\documentclass{amsart}
\usepackage{amsmath}
\usepackage{amsthm}
\usepackage{amsfonts}
\usepackage{amssymb}
\usepackage{hyperref,color}
\usepackage{array}
\usepackage{empheq}
\usepackage{tikz}
\usetikzlibrary{matrix,arrows,positioning,automata,shapes}
\usepackage{tikz-qtree}
\tikzset{
   level distance=1.2cm,sibling distance=.8cm, 
   edge from parent path={(\tikzparentnode) -- (\tikzchildnode)},
   elb/.style={draw,ellipse,inner sep=1pt}, 
   wlb/.style={fill=white,inner sep=1pt}, 
   glb/.style={fill=black!7,inner sep=1pt}, 
   gelb/.style={draw,ellipse,fill=black!7,inner sep=1pt}, 
   nlb/.style={inner sep=0pt,draw=none,minimum size=0pt}, 
   highlight/.style={line width=.5cm,color=black!22,cap=round,join=round,opacity=0.3}
}
\pgfdeclarelayer{background layer}
\pgfsetlayers{background layer,main}
\definecolor{bkg}{rgb}{0.8,0.8,0.8}

\newcommand{\mB}[1]{\mathbf{#1}}
\newcommand{\Bwr}{\boldsymbol{\wr}}
\newcommand{\PF}{\noindent {\it Proof:\ }}
\newcommand{\QED}{\hfill $\Box$}

\newtheorem{theorem}{Theorem}[section]
\theoremstyle{definition}
\newtheorem{definition}[theorem]{Definition}
\newtheorem{example}{Example}

\newtheorem{fact}{Fact}

\newcommand{\Z}{\mathbb{Z}}

\newcommand{\lp}{\left(}
\newcommand{\rp}{\right)}

\newcommand{\cC}{\mathcal{C}}

\DeclareMathOperator{\Aut}{\text{Aut}}

\newcommand{\SgpDec}{\textsc{SgpDec}}
\newcommand{\GAP}{\textsc{Gap}}

\title{Cascade Product of Permutation Groups}
\author{Attila Egri-Nagy$^{1,2}$ and Chrystopher L. Nehaniv$^1$}
\address{$^1$Centre for Computer Science \& Informatics Research,
        University of Hertfordshire, College Lane, Hatfield, Herts AL10 9AB, United Kingdom
\and
$^2$Centre for Research in Mathematics, School of Computing, Engineering and Mathematics, University of Western Sydney (Parramatta Campus), Locked Bag 1797, Penrith, NSW 2751, Australia}
\email{A.Egri-Nagy@uws.edu.au,\ C.L.Nehaniv@herts.ac.uk}

\begin{document}
\maketitle

\begin{abstract}
Motivated by computational efficiency in algebraic automata theory here we define the cascade product of permutation groups as an external product, as a generic extension.
It is the most general hierarchical product that uses arbitrary total functions
to combine a linearly ordered set of permutation groups.
Algebraically speaking, cascade products are explicitly constructed substructures of the iterated wreath product.
We show how direct, semidirect and wreath products can be described as cascade products and we also discuss examples of composite groups that can only be constructed exactly as generic extensions by cascade products.
The cascade construction naturally generalizes to the transformation semigroup
case by leaving out the details of defining inverse operations.
\end{abstract}


\section{Introduction}

The idea, that connections between component groups in a product can be represented by functions, implicitly appears in semidirect and wreath products.
In mathematics, there was no apparent reason to pay attention to these functions.
In computer science, however, we need data structures to describe the internal structure of these compositions.
We need to specify what to do in a component based on what happened in another one
by using functions, the so-called dependency functions.
With these functions we can define cascade products simply by restricting the dependencies to be unidirectional.
In turn, the cascade product allows us to construct the usual products and more generic extensions, based on what sort of functions we allow as dependencies.
So, at the end of the day, the dependency functions are actually mathematically interesting objects.

The idea of dependency functions is implicit in the holonomy decomposition of transformation semigroups (see e.g.\ \cite{automatanetworks2005}) and sporadically appeared in some disguise in automata theory \cite{Nozaki78} and in statistics for ordering experiments~\cite{GeneralizedWreathProducts1983}, but the arrival of computational implementations of hierarchical semigroup decompositions made it indispensable. Here we concentrate only on groups. The semigroup case is easier in a sense that  we do not need to characterize inverses.

The main purpose of this tutorial is to give mathematically precise
documentation for the computational implementation of cascade products in
\SgpDec \cite{sgpdec}, which is a \GAP \cite{GAP4} package for hierarchical
decompositions of permutation groups and transformation semigroups.
In order to make the text more readable we discuss interesting examples.
We also think that this explicit construction of cascade products yields new questions and will eventually lead to new mathematical results. 
\subsection{Notation}
We denote the identity of a group $G$ by $1$ or use $1_G$ if there is ambiguity.
Let $X$ be a set, called the \emph{state set}, and  a group $G$ acting on this
set, meaning that $x^1=x$ and $\lp x^g\rp^h=x^{gh}$ for all $x\in X$ and $g,h\in
G$. We denote this permutation group by $(X,G)$. We call $|X|$ the
\emph{degree}, and $|G|$ the \emph{order} of the permutation group.
The \emph{order of a group element} $g$ is the size of the group it generates,
$|\langle g \rangle|$.
This is also called as the \emph{period length}, as in finite groups the powers of a
single element form a cycle.

Let $\mB{n}$ be the set $\{1,\ldots,n\}$. A mod-$n$ \emph{counter} is $(\mB{n},
\Z_n)$, the cyclic group $\Z_n=\langle +1\rangle$ acting on $\mB{n}$. For the
counters only, to represent addition modulo $n$, it is useful to relabel $\mB{n}$ to $\{0,1,\ldots,n-1\}$. By $S_n$ we denote the \emph{symmetric group}, the group of all permutations of $n$ points. We use $D_n$ for the \emph{dihedral group}, the symmetry group of the regular $n$-gon.
For more on permutation groups see \cite{DixonMortimerPermGroups96,CameronPermGroups99}.

\subsection{Motivating Example}
We would like to build the mod-4 counter from two mod-2 counters. The direct
product $\Z_2\times \Z_2$ contains no element of order 4, therefore it cannot be
used as a mod-4 counter.
Since the automorphism group Aut$(\Z_2)$ is trivial, there is only one semidirect product, which reverts back to the direct product. If we try the wreath product, then we get $\Z_2\wr \Z_2 \cong D_4$,
the dihedral group of the square, which is twice the size of $\Z_4$. The embedding obviously works $\Z_4\hookrightarrow D_4$, thus the wreath product can be used as a mod-4 counter. But from the engineering perspective, this construction is not efficient. The dihedral group also has the flip-symmetry, in addition to the required rotations. 
Therefore we need a product that allows us to build $(\mB{4},\Z_4)$ exactly, i.e.\ as an isomorphic composite structure.

\section{Cascade Products of Permutation Groups}

Let $L=\left[(X_1,G_1),\ldots,(X_n,G_n)\right]$ be an ordered list of
permutation groups.
These are called the \emph{components}, and we build another permutation group
from them.
In general, we can build many different cascade products using the same list of components. 

The components  in a cascade product are connected according to the hierarchy defined by the given linear order of the indices $1<2<..<n$.
We say that $(X_j,G_j)$ \emph{depends on}, is connected to, or influenced by
$(X_i,G_i)$ if $i<j$.
Note that this formal dependence does not imply actual control information flow, as
in the case of direct product.
The index encodes the \emph{depth} of the component\footnote{This left-to-right top-down ordering is due to the constraints of a software implementation, as in computer algebra systems lists are usually indexed by starting from 1 and most algorithms dealing with cascade products start from the independent component.}.
Thus, $(X_1,G_1)$ is  the \emph{top} level, the independent component, while $(X_n,G_n)$ is the most dependent \emph{bottom} level component.

\subsection{Coordinates: States of  a Cascade Product}
The set of states on which a cascade product acts on is $X_1\times\cdots\times
X_n$, thus a point or a state is an $n$-tuple $x=(x_1,\dots,x_n)$, $x_i\in X_i$.
We call the values in $x$  the \emph{coordinates} and say that $x_i$ is the
$i$th level coordinate. 

\subsection{Dependency Functions}
The elements of a cascade product will act as permutations of $X_1\times\cdots\times X_n$, but not just any permutations. Only those are allowed that respect the hierarchical order of the components. Mathematically it would be enough to define the set of allowed permutations by this property, saying that the component actions lie within the components. However, computationally we have to construct the operations explicitly.
In order to ensure that the permutations in a cascade product respect the hierarchical order, we describe them in terms of functions producing elements of  component groups based on states in other components above.

\begin{definition}
A \emph{dependency function} $d_i$ of level $i$ for a list of components $L$ is a function
$$d_i: X_1\times\cdots\times X_{i-1}\rightarrow G_i,\ \ \ i\in \mB{n}.$$
\end{definition}
\noindent A dependency function of level $i$ takes $i-1$ arguments. On the first level we have the empty product therefore we define $(x_1,\ldots,x_{i-1})=\varnothing$ for $i=1$. If no confusion arises, on the top level we can simply write $d_1\in G_1$ instead of $d_1(\varnothing)\in G_1$. 
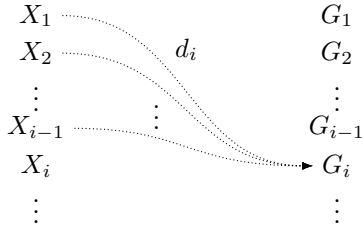
\begin{center}\tikzset{>=latex}
\begin{tikzpicture}
\newcommand\sh{4} 
\newcommand\x{.5}  
\draw 
 (0,4*\x) node (X1) {$X_1$}
 (0,3*\x) node (X2) {$X_2$}
 (0,2*\x) node {$\vdots$}
 (0,\x) node (Xi-1) {$X_{i-1}$}
 (0,0) node (Xi) {$X_i$}
 (0,-\x) node {$\vdots$}
 (\sh,4*\x) node (S1) {$G_1$}
 (\sh,3*\x) node (S2) {$G_2$}
 (\sh,2*\x) node {$\vdots$}
 (\sh,\x) node (Si-1) {$G_{i-1}$}
 (\sh,0) node (Si) {$G_i$}
 (\sh,-\x) node {$\vdots$};
\draw[densely dotted,->] (X1)..controls (.5*\sh,4*\x) and (.5*\sh,0*\x) ..(Si);
\draw[densely dotted,->] (X2)..controls (.5*\sh,3*\x) and (.5*\sh,0*\x) ..(Si);
\draw[densely dotted,->] (Xi-1)..controls (.5*\sh,\x) and (.5*\sh,0*\x) ..(Si);
\draw (.4*\sh,1.5*\x) node {$\vdots$};
\draw (.5*\sh,3.1*\x) node {$d_i$};
\end{tikzpicture}\end{center}
The intuitive interpretation of the dependency function is that it is a lookup table from which one can pick an actual permutation from a component based on the states of the components above. We also say that a dependency function encodes the connections into one particular component.

\subsection{Permutation Cascades}
\begin{definition}
A \emph{permutation cascade} for a given list of components $L$ is an $n$-tuple of dependency functions $(d_1,\ldots,d_n)$, where $d_i$ is a dependency function of level $i$.
\end{definition}
Cascades can be visualized as labelled trees: the branches encode the dependency arguments and the nodes contain the values, the local actions.
\begin{example}
A cascade  $c=(d_1,d_2,d_3)$ acting on $\mB{2}\times\mB{2}\times\mB{3}$ can be drawn as the following labelled tree.
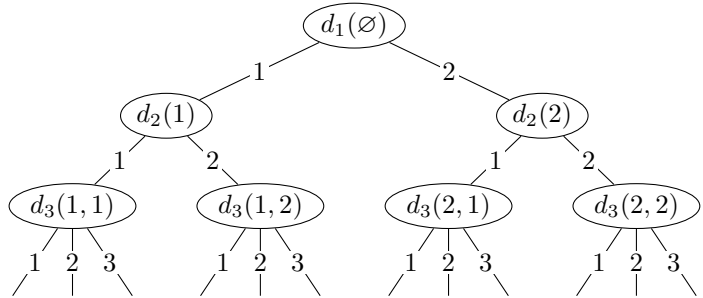
\begin{center}
\begin{tikzpicture}
[every tree node/.style=nlb]
\Tree [.\node[elb] {$d_1(\varnothing)$};
  \edge node[wlb] {1}; [.\node[elb]  {$d_2(1)$};
    \edge node[wlb] {1}; [.\node[elb] {$d_3(1,1)$}; 
       \edge node[wlb] {1};[.{} ]
       \edge node[wlb] {2};[.{} ]
       \edge node[wlb] {3};[.{} ] ] 
     \edge node[wlb] {2};[.\node[elb]  {$d_3(1,2)$};
       \edge node[wlb] {1};[.{} ]
       \edge node[wlb] {2};[.{} ]
       \edge node[wlb] {3}; [.{} ] ] ]
   \edge node[wlb] {2}; [.\node[elb] {$d_2(2)$};
     \edge node[wlb] {1}; [.\node[elb] {$d_3(2,1)$};
       \edge node[wlb] {1};[.{} ]
       \edge node[wlb] {2};[.{} ]
       \edge node[wlb] {3};[.{} ] ] 
     \edge node[wlb] {2};[.\node[elb] {$d_3(2,2)$};
       \edge node[wlb] {1};[.{} ]
       \edge node[wlb] {2};[.{} ]
       \edge node[wlb] {3};[.{} ] ] ]
]
\end{tikzpicture}
\end{center}
\end{example}
\subsection{Coordinatewise Action} The action of a permutation cascade  $d=(d_1, \ldots, d_n)$  on coordinates $x=(x_1,\ldots, x_n)$ for level $i$ is defined  by 
\begin{align*}
 x_i^d&:=x_i^{d_i(x_1,\ldots,x_{i-1})}
\end{align*}
thus the full action is
\begin{equation*}
x^d=(x_1,\ldots,x_n)^{(d_1,\ldots,d_n)}=\left(x_1^{d_1(\varnothing)},x_2^{d_2(x_1)},\ldots,x_n^{d_n(x_1,\ldots,x_{n-1})}\right).
\end{equation*}

\begin{example}[Cascade action]
\label{ex:3levelpermcascadeaction}
Let's consider the cascade defined by the following dependencies.
In order to avoid confusion, we denote coordinate tuples by square brackets.
$\varnothing\mapsto(1,2)$,

$[1]\mapsto()$, 
$[2]\mapsto(1,2)$, 

$[1,1]\mapsto(1,2,3)$, 
$[1,2]\mapsto(2,3)$, 
$[2,1]\mapsto(1,3)$, 
$[2,2]\mapsto(1,3,2)$ 

\begin{center}
\begin{tikzpicture}
[every tree node/.style=nlb]
\Tree [.\node[elb] (b) {(1,2)};
  \edge node[wlb] {1}; [.\node[elb] (a)  {()};
    \edge node[wlb] {1}; [.\node[elb] (d) {(1,2,3)}; 
       \edge node[wlb] {1};[.{} ]
       \edge node[wlb] {2};[.{} ]
       \edge node[wlb] {3};[.{} ] ] 
     \edge node[wlb] {2};[.\node[elb] (c) {(2,3)};
       \edge node[wlb] {1};[.{} ]
       \edge node[wlb] {2};[.{} ]
       \edge node[wlb] {3}; [.{} ] ] ]
   \edge node[wlb] {2}; [.\node[elb] {(1,2)};
     \edge node[wlb] {1}; [.\node[elb] {(1,3)};
       \edge node[wlb] {1};[.{} ]
       \edge node[wlb] {2};[.{} ]
       \edge node[wlb] {3};[.{} ] ] 
     \edge node[wlb] {2};[.\node[elb] {(1,3,2)};
       \edge node[wlb] {1};[.{} ]
       \edge node[wlb] {2};[.{} ]
       \edge node[wlb] {3};[.{} ] ] ]
]
\end{tikzpicture}
\end{center}
In order to evaluate the action of this cascade we need to carry out the local in action each node corresponding to subtrees.
We do only relabelling on the local branches restricted a single level (as opposed to moving subtrees recursively).
Theoretically the local actions are done synchronously, but we can still do them sequentially. 
Fortunately, the order does not matter. 
For instance, we can do it top-down, starting at the first level,
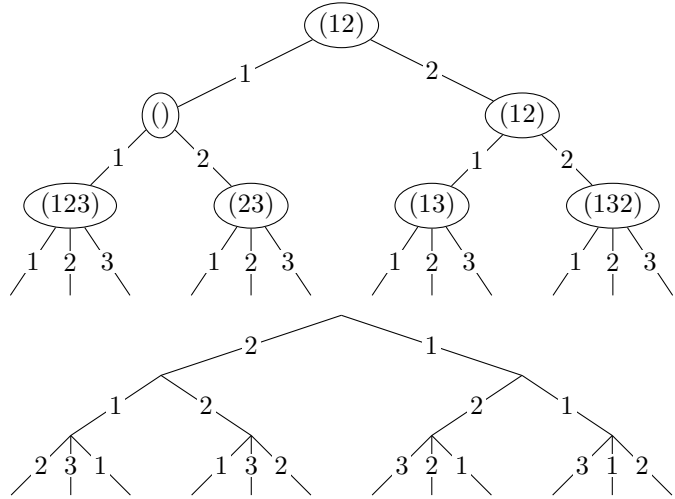
\begin{center}
\begin{tikzpicture}
[every tree node/.style=nlb]
\Tree [
  \edge node[wlb] {2}; [.\node[elb] (a)  {()};
    \edge node[wlb] {1}; [.\node[elb] (d) {(1,2,3)}; 
       \edge node[wlb] {1};[.{} ]
       \edge node[wlb] {2};[.{} ]
       \edge node[wlb] {3};[.{} ] ] 
     \edge node[wlb] {2};[.\node[elb] (c) {(2,3)};
       \edge node[wlb] {1};[.{} ]
       \edge node[wlb] {2};[.{} ]
       \edge node[wlb] {3}; [.{} ] ] ]
   \edge node[wlb] {1}; [.\node[elb] {(1,2)};
     \edge node[wlb] {1}; [.\node[elb] {(1,3)};
       \edge node[wlb] {1};[.{} ]
       \edge node[wlb] {2};[.{} ]
       \edge node[wlb] {3};[.{} ] ] 
     \edge node[wlb] {2};[.\node[elb] {(1,3,2)};
       \edge node[wlb] {1};[.{} ]
       \edge node[wlb] {2};[.{} ]
       \edge node[wlb] {3};[.{} ] ] ]
]
\end{tikzpicture}
\end{center}
then proceeding to the second level,
\begin{center}
\begin{tikzpicture}
[every tree node/.style=nlb]
\Tree [
  \edge node[wlb] {2}; [
    \edge node[wlb] {1}; [.\node[elb]  {(1,2,3)}; 
       \edge node[wlb] {1};[.{} ]
       \edge node[wlb] {2};[.{} ]
       \edge node[wlb] {3};[.{} ] ] 
     \edge node[wlb] {2};[.\node[elb]  {(2,3)};
       \edge node[wlb] {1};[.{} ]
       \edge node[wlb] {2};[.{} ]
       \edge node[wlb] {3}; [.{} ] ] ]
   \edge node[wlb] {1}; [
     \edge node[wlb] {2}; [.\node[elb] {(1,3)};
       \edge node[wlb] {1};[.{} ]
       \edge node[wlb] {2};[.{} ]
       \edge node[wlb] {3};[.{} ] ] 
     \edge node[wlb] {1};[.\node[elb] {(1,3,2)};
       \edge node[wlb] {1};[.{} ]
       \edge node[wlb] {2};[.{} ]
       \edge node[wlb] {3};[.{} ] ] ]
]
\end{tikzpicture}
\end{center}
finally we do the relabellings on the lowest levels.
\begin{center}
\begin{tikzpicture}
[every tree node/.style=nlb]
\Tree [
  \edge node[wlb] {2}; [
    \edge node[wlb] {1}; [
       \edge node[wlb] {2};[.{} ]
       \edge node[wlb] {3};[.{} ]
       \edge node[wlb] {1};[.{} ] ] 
     \edge node[wlb] {2};[
       \edge node[wlb] {1};[.{} ]
       \edge node[wlb] {3};[.{} ]
       \edge node[wlb] {2}; [.{} ] ] ]
   \edge node[wlb] {1}; [
     \edge node[wlb] {2}; [
       \edge node[wlb] {3};[.{} ]
       \edge node[wlb] {2};[.{} ]
       \edge node[wlb] {1};[.{} ] ] 
     \edge node[wlb] {1};[
       \edge node[wlb] {3};[.{} ]
       \edge node[wlb] {1};[.{} ]
       \edge node[wlb] {2};[.{} ] ] ]
]
\end{tikzpicture}
\end{center}
\end{example}

\subsection{Multiplication}

Let $d=(d_1,\ldots,d_n)$ and $f=(f_1,\ldots,f_n)$ be two permutation cascades defined over the same list of components. Multiplication goes along the combination of functions but we have to keep track of the moving arguments. We combine on the right. The coordinatewise action of $df$ on $x=(x_1,\ldots,x_n)$ is
$$\lp x_i^d\rp^f=\left(x_i^{d_i(x_1,\ldots,x_{i-1})}\right)^{f_i\lp x_1^d,\ldots,x_{i-1}^d\rp}=x_i^{d_i(x_1,\ldots,x_{i-1})f_i\lp x_1^d,\ldots,x_{i-1}^d\rp}=x_i^{df}$$
and from the exponent we can express the combined dependency function:
$$(df)_i:(x_1,\ldots,x_{i-1})\mapsto d_i(x_1,\ldots,x_{i-1})f_i\lp x_1^d,\ldots,x_{i-1}^d\rp $$

Associativity follows from the associativity of  the group multiplication and of the group action. If $h=(h_1,\ldots,h_{n})$ is another permutation cascade, then
\begin{align*}
\left((df)h\right)_i &= \left(d_i(x_1,\ldots,x_{i-1})f_i\lp x_1^d,\ldots,x_{i-1}^d\rp\right)h_i\lp x_1^{df},\ldots,x_{i-1}^{df}\rp\\
 &= d_i(x_1,\ldots,x_{i-1})f_i\lp x_1^d,\ldots,x_{i-1}^d\rp h_i\lp x_1^{df},\ldots,x_{i-1}^{df}\rp\\
&= d_i(x_1,\ldots,x_{i-1})\left(f_i\lp x_1^d,\ldots,x_{i-1}^d\rp h_i\lp\lp x_1^d\rp^f,\ldots,\lp x_{i-1}^d\rp^f\rp\right)\\
&=\left(d(fh)\right)_i.
\end{align*}

\begin{example}[Multiplying permutation cascades in tree form]
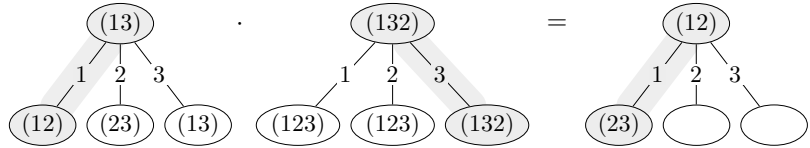
\label{ex:cascmul}

\begin{center}
\begin{tikzpicture}
[level distance=2cm,sibling distance=.15cm]
\tikzstyle{label} =[fill=white,inner sep=3pt]
\Tree [.\node[gelb] (r) {(1,3)};
        \edge node[glb] {1}; [.\node[gelb] (r1)  {(1,2)}; ]
        \edge node[wlb] {2}; [.\node[elb]  {(2,3)}; ]
        \edge node[wlb] {3}; [.\node[elb]  {(1,3)}; ]
]
\begin{pgfonlayer}{background layer}
\draw[highlight] (r)--(r1);
\end{pgfonlayer}
\end{tikzpicture}
\begin{tikzpicture}
\node (a) at (0,1.5) {$\cdot$};
\node (n) at (0,0) {};
\end{tikzpicture}
\begin{tikzpicture}
[level distance=2cm,sibling distance=.15cm]
\Tree [.\node[gelb] (r) {(1,3,2)};
        \edge node[wlb] {1}; [.\node[elb]   {(1,2,3)}; ]
        \edge node[wlb] {2}; [.\node[elb]  {(1,2,3)}; ]
        \edge node[glb] {3}; [.\node[gelb] (r3) {(1,3,2)}; ]
]
\begin{pgfonlayer}{background layer}
\draw[highlight] (r)--(r3);
\end{pgfonlayer}
\end{tikzpicture}
\begin{tikzpicture}
\node (a) at (0,1.5) {=};
\node (n) at (0,0) {};
\end{tikzpicture}
\vskip3pt
\begin{tikzpicture}
[level distance=2cm,sibling distance=.15cm]
\tikzstyle{label} =[fill=white,inner sep=3pt]
\Tree [.\node[gelb] (r) {(1,2)};
        \edge node[glb] {1}; [.\node[gelb] (r1)  {(2,3)}; ]
        \edge node[wlb] {2}; [.\node[elb]  {\color{white}(1,2)}; ]
        \edge node[wlb] {3}; [.\node[elb] {\color{white}(2,3)}; ]
]
\begin{pgfonlayer}{background layer}
\draw[highlight] (r)--(r1);
\end{pgfonlayer}
\end{tikzpicture}
\end{center}
Calculating the highlighted path in the product, first we multiply the top level permutations: $(1,3)\cdot(1,3,2)=(1,2)$.
Then we take the level 2 permutation along the path $[1]$ from the first, and the level 2 permutation along path $[3]$ from the second cascade, since $[1]$ is moved to $[3]$ by the top level permutation of the first cascade. 
Thus $(1,2)\cdot(1,3,2)=(2,3)$.

Then, path [2] is fixed by $c_1$, so in this particular case we simply just multiply the corresponding local actions.
Note that we don't need calculate the top level again, only $(2,3)\cdot(1,2,3)=(1,2)$. 
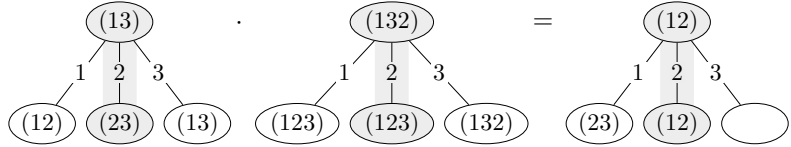
\begin{center}
\begin{tikzpicture}
[level distance=2cm,sibling distance=.15cm]
\tikzstyle{label} =[fill=white,inner sep=3pt]
\Tree [.\node[gelb] (r) {(1,3)};
        \edge node[wlb] {1}; [.\node[elb] (r1)  {(1,2)}; ]
        \edge node[glb] {2}; [.\node[gelb] (r2) {(2,3)}; ]
        \edge node[wlb] {3}; [.\node[elb]  {(1,3)}; ]
]
\begin{pgfonlayer}{background layer}
\draw[highlight] (r)--(r2);
\end{pgfonlayer}
\end{tikzpicture}
\begin{tikzpicture}
\node (a) at (0,1.5) {$\cdot$};
\node (n) at (0,0) {};
\end{tikzpicture}
\begin{tikzpicture}
[level distance=2cm,sibling distance=.15cm]
\Tree [.\node[gelb] (r) {(1,3,2)};
        \edge node[wlb] {1}; [.\node[elb]   {(1,2,3)}; ]
        \edge node[glb] {2}; [.\node[gelb] (r2) {(1,2,3)}; ]
        \edge node[wlb] {3}; [.\node[elb]  {(1,3,2)}; ]
]
\begin{pgfonlayer}{background layer}
\draw[highlight] (r)--(r2);
\end{pgfonlayer}
\end{tikzpicture}
\begin{tikzpicture}
\node (a) at (0,1.5) {=};
\node (n) at (0,0) {};
\end{tikzpicture}
\vskip3pt
\begin{tikzpicture}
[level distance=2cm,sibling distance=.15cm]
\tikzstyle{label} =[fill=white,inner sep=3pt]
\Tree [.\node[gelb] (r) {(1,2)};
        \edge node[wlb] {1}; [.\node[elb]   {(2,3)}; ]
        \edge node[glb] {2}; [.\node[gelb] (r2)  {(1,2)}; ]
        \edge node[wlb] {3}; [.\node[elb] {\color{white}(2,3)}; ]
]
\begin{pgfonlayer}{background layer}
\draw[highlight] (r)--(r2);
\end{pgfonlayer}
\end{tikzpicture}
\end{center}

Finally, path [3] is sent to [1], so $(1,3)\cdot(1,2,3)=(2,3)$ finishes the calculation.
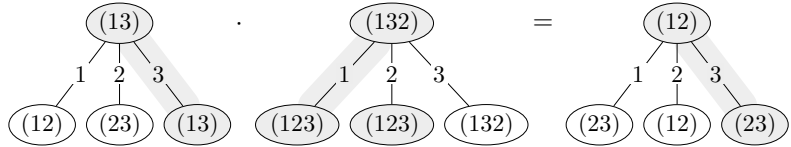
\begin{center}
\begin{tikzpicture}
[level distance=2cm,sibling distance=.15cm]
\tikzstyle{label} =[fill=white,inner sep=3pt]
\Tree [.\node[gelb] (r) {(1,3)};
        \edge node[wlb] {1}; [.\node[elb] (r1)  {(1,2)}; ]
        \edge node[wlb] {2}; [.\node[elb] (r2) {(2,3)}; ]
        \edge node[glb] {3}; [.\node[gelb] (r3) {(1,3)}; ]
]
\begin{pgfonlayer}{background layer}
\draw[highlight] (r)--(r3);
\end{pgfonlayer}
\end{tikzpicture}
\begin{tikzpicture}
\node (a) at (0,1.5) {$\cdot$};
\node (n) at (0,0) {};
\end{tikzpicture}
\begin{tikzpicture}
[level distance=2cm,sibling distance=.15cm]
\Tree [.\node[gelb] (r) {(1,3,2)};
        \edge node[glb] {1}; [.\node[gelb] (r1)  {(1,2,3)}; ]
        \edge node[wlb] {2}; [.\node[gelb] (r2) {(1,2,3)}; ]
        \edge node[wlb] {3}; [.\node[elb]  {(1,3,2)}; ]
]
\begin{pgfonlayer}{background layer}
\draw[highlight] (r)--(r1);
\end{pgfonlayer}
\end{tikzpicture}
\begin{tikzpicture}
\node (a) at (0,1.5) {=};
\node (n) at (0,0) {};
\end{tikzpicture}
\vskip3pt
\begin{tikzpicture}
[level distance=2cm,sibling distance=.15cm]
\tikzstyle{label} =[fill=white,inner sep=3pt]
\Tree [.\node[gelb] (r) {(1,2)};
        \edge node[wlb] {1}; [.\node[elb]   {(2,3)}; ]
        \edge node[wlb] {2}; [.\node[elb] (r2)  {(1,2)}; ]
        \edge node[glb] {3}; [.\node[gelb] (r3){(2,3)}; ]
]
\begin{pgfonlayer}{background layer}
\draw[highlight] (r)--(r3);
\end{pgfonlayer}
\end{tikzpicture}
\end{center}
\end{example}

\subsection{Identity} The identity of a cascade product is $1=(e_1,e_2,\ldots, e_n)$ where $e_i$ is a constant function with $e_i:(x_1, \ldots, x_{i-1})\mapsto 1_{G_i}$ for all $(x_1, \ldots, x_{i-1})\in X_1\times \ldots \times X_{i-1}$.
\subsection{Inverses} For finding the left inverse of a permutation cascade $d$ we use
$$x_i^{fd}=x_i^{f_i(x_1,\ldots,x_{i-1})d_i\lp x_1^f,\ldots,x_{i-1}^f\rp}$$
so in order to make $f=d^{-1}$ we should define $f$ by
$$f_i(x_1,\ldots,x_{i-1}):=\left(d_i\lp x_1^f,\ldots,x_{i-1}^f\rp\right)^{-1}$$
which first looks like a circular definition but it is actually a recursive one. By explicitly writing down the dependency functions we get
\begin{align*}
f_1 &:= d_1(\varnothing)^{-1}\\
f_2(x_1) &:= \left(d_2(x_1^{f_1})\right)^{-1}=\left(d_2(x_1^{d_1(\varnothing)^{-1}})\right)^{-1}\\
f_3(x_1,x_2) &:= \left(d_3\left(x_1^{f_1},x_2^{f_2(x_1)}\right) \right)^{-1}=\left(d_3\left(x_1^{d_1(\varnothing)^{-1}},x_2^{\left(d_2\left(x_1^{d_1(\varnothing)^{-1}}\right)\right)^{-1}}\right) \right)^{-1}\\
&\vdots
\end{align*}
therefore $f$ is a left inverse. Now we show that it is also a right inverse. Multiplying $d$ by $f$ on the right
\begin{equation*}
x_i^{df}=x_i^{d_i(x_1,\ldots,x_{i-1})f_i\lp x_1^d,\ldots,x_{i-1}^d\rp}
\end{equation*}
we only have to evaluate $f_i\left(x_1^d,\ldots,x_{i-1}^d\right)$. Doing it recursively:
\begin{align*}
f_1 &= d_1(\varnothing)^{-1}\\
f_2(x_1^d) &=\left(d_2\left(\left(x_1^{d_1(\varnothing)}\right)^{d_1(\varnothing)^{-1}}\right)\right)^{-1}=\left(d_2\left(x_1^{d_1(\varnothing)d_1(\varnothing)^{-1}}\right)\right)^{-1}=\lp d_2\lp x_1\rp\rp^{-1}\\
f_3(x_1^d,x_2^d) &= \lp d_3\lp \lp x_1^{d_1(\varnothing)}\rp^{f_1},\lp x_2^{d_2(x_1)}\rp^{f_2(x_1)}\rp \rp^{-1}\\
&= \lp d_3\lp \lp x_1^{d_1(\varnothing)}\rp^{d_1(\varnothing)^{-1}},\lp x_2^{d_2(x_1)}\rp^{\lp d_2\lp x_1\rp\rp^{-1}}\rp \rp^{-1}\\
&= \lp d_3\lp x_1^{d_1(\varnothing)d_1(\varnothing)^{-1}}, x_2^{d_2(x_1)\lp d_2\lp x_1\rp\rp^{-1}}\rp \rp^{-1}= \lp d_3\lp x_1, x_2\rp \rp^{-1}\\
&\vdots
\end{align*}
The pattern can easily be seen and by induction we can show that it is true for each level. Supposing that $f_j\lp x_1^d,\ldots,x_{j-1}^d\rp=\lp d_j(x_1^f,\ldots,x_{j-1}^f)\rp^{-1}$ for all $j$ from 1 to $i$, we have $x_j^{df}=x_j$. Then
$$f_{i+1} \lp x_1^d,\ldots,x_{i}^d \rp = \lp d_{i+1} \lp x_1^{df},\ldots,x_{i}^{df}\rp  \rp ^{-1} = \lp d_{i+1} \lp x_1,\ldots,x_{i}\rp  \rp ^{-1}.$$
It follows immediately that $f$ is a right-, hence unique two-sided inverse of $d$.
We have proved the following
\begin{theorem}[Computing Inverses of Permutation Cascades]
If $d=(d_1,\ldots,d_n)$ is a permutation cascade, then its inverse is $f=(f_1,\ldots,f_n)$, where each
$$f_i(x_1^d,\ldots,x_{i-1}^d):=\left(d_i(x_1,\ldots,x_{i-1})\right)^{-1}.$$
\end{theorem}
\subsection{The Full Cascade Product}
Let $\cC_L$ be a set of all permutation cascades defined on the list of components $L=\left[(X_1,G_1),\ldots,(X_n,G_n)\right]$. Then the  \emph{full cascade product} defined on $L$  is the permutation group
$$(X_1,G_1) \Bwr\cdots\Bwr (X_n,G_n) := (X_1\times\cdots\times X_n, \cC_L).$$
\noindent The degree of the full cascade product is $\prod_{i=1}^n|X_i|$ and the order can be calculated by
$$|\cC_L|=|G_1|\prod_{i=2}^n|G_i|^{\prod_{j=1}^{i-1}|X_j|}. $$
\noindent Let $o_i$ be the order and $n_i$ be now the degree of $(X_i,G_i)$, then
$$|(X_1,G_1)\Bwr(X_2,G_2)\Bwr(X_3,G_3)|=o_1\cdot o_2^{n_1}\cdot o_3^{n_1n_2}.$$
Let $L_1=\left[(X_1,G_1),(X_2,G_2)\right]$ and $L_2=\left[(X_2,G_2),\ldots,(X_3,G_3)\right]$, then
\begin{align*}
|\cC_{L_1}| &= o_1\cdot o_2^{n_1}\\
|\cC_{L_2}| &= o_2\cdot o_3^{n_2}
\end{align*}
Therefore,
\begin{align*}
|(X_1\times X_2,\cC_{L_1})\Bwr(X_3,G_3)|&=o_1o_2^{n_1}\cdot o_3^{n_1n_2}\\
&=o_1\cdot(o_2o_3^{n_2})^{n_1}\\
&=|(X_1,G_1)\Bwr(X_2\times X_3,\cC_{L_2})|.
\end{align*}
\begin{fact}
The full cascade product is associative.
\begin{align*}
\left((X_1,G_1)\Bwr(X_2,G_2)\right)\Bwr(X_3,G_3)&\cong(X_1,G_1)\Bwr((X_2,G_2)\Bwr(X_3,G_3))\\
&\cong (X_1,G_1)\Bwr(X_2,G_2)\Bwr(X_3,G_3).
\end{align*}
\end{fact}
The full cascade product behaves similarly as the iterated permutation wreath product, actually it is isomorphic to that (see section \ref{sec:wreath}). Therefore, we use symbols $\Bwr$ for the full cascade and $\wr$ for the wreath product. 

\begin{example}
$(\mB{3},S_3)\Bwr(\mB{2},\Z_2)\cong \Z_2\times S_4$, a group of order 48. However, $(S_3,S_3)\Bwr(\mB{2},\Z_2)$ is a group of size 384. This shows that isomorphic component groups may yield different cascade products if they act on different sets. 
\end{example}

\subsection{Cascade Products} 
\begin{figure}
\begin{center}
\tikzset{->,>=triangle 45,auto,node distance=4cm}
\begin{tikzpicture}
\tikzstyle{every state}=[minimum size=3pt]
  \node[state]  (q_1)  {$0$};
  \node[fill=bkg]  (c0) [below left=0.7cm and 0.8cm of q_1] {$c(0)$};
  \node[state]  (q_2) [right of=q_1] {$1$};
  \node () [right=0.5cm of q_2] {top};
  \node[fill=bkg]  (c1) [below left=0.7cm and 0.5cm of q_2] {$c(1)$};
  \node[state]  (q_3) [below of=q_1] {$0$};
  \node[state]  (q_4) [right of=q_3] {$1$};
  \node () [right=0.5cm of q_4] {bottom};

  \path (q_1) edge  [bend left] node {$+1$} (q_2);
  \path (q_1) edge  [loop] node  [above] {$1_{\Z_2}$} (q_1);
  \path (q_2) edge  [bend left] node [above] {$+1$} (q_1);
  \path (q_2) edge  [loop] node  [above] {$1_{\Z_2}$} (q_2);
  \path (q_3) edge  [loop] node (id3) [above] {$1_{\Z_2}$} (q_3);
  \path (q_4) edge  [loop] node (id4) [above] {$1_{\Z_2}$} (q_4);
  \path (q_3) edge  [bend left] node (ct) {$+1$} (q_4);
  \path (q_4) edge  [bend left] node (cb) [above] {$+1$} (q_3);
  \path[>=angle 60,densely dotted]  (q_1) edge (c0);
  \path[>=angle 60,densely dotted]  (c0) edge (id3);
  \path[>=angle 60,densely dotted]  (c0) edge (id4);
  \path[>=angle 60,densely dotted]  (q_2) edge (c1);
  \path[>=angle 60,densely dotted]  (c1) edge (ct);
  \path[>=angle 60,densely dotted]  (c1) edge (cb);
\end{tikzpicture}
\caption{Two mod-2 counters coupled together to build a mod-4 counter.}
\label{fig:mod4counter}
\end{center}
\end{figure}
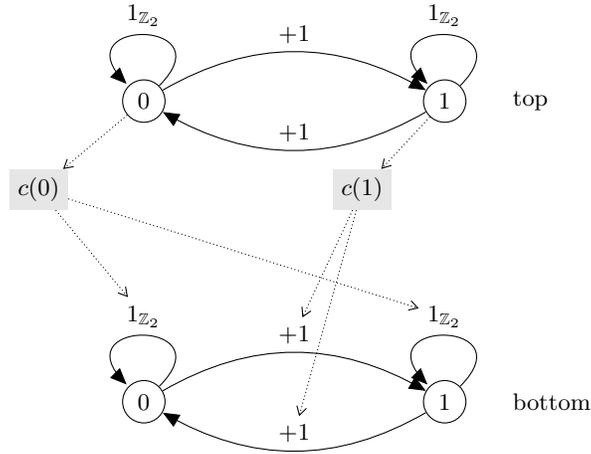

Let $L$ be a given list of components $\left[(X_1,G_1),\ldots,(X_n,G_n)\right]$. Then a \emph{cascade product} is a permutation group
$$(X_1\times\cdots\times X_n, G)$$
\noindent where $G\leq\cC_L$, i.e. a subgroup of the full cascade product. The most usual scenario is that we have a generating set $W$ of permutation cascades, $G=\langle W\rangle$. In this case we can use the notation:
$$(X_1,G_1) \wr_W\cdots\wr_W (X_n,G_n) := (X_1\times\cdots\times X_n, \langle W\rangle).$$
The meaning of the $\wr_W$ symbol is that the connections are restricted to what is generated by $W$.
\begin{example} We can construct $(\mB{4},\Z_4)$ exactly by using two copies of $(\mB{2},\Z_2)$. The generator set contains only one permutation cascade  $W=\{(+1,c)\}$, where $+1$ is the generator of $\Z_2$ and $c$ is a dependency function with mappings 
\begin{align*}
0&\mapsto 1_{\Z_2},\\
1&\mapsto +1.
\end{align*}
The first dependency is trivial while the second dependency implements the carry.
$$(\mB{2},\Z_2)\wr_W(\mB{2},\Z_2)\cong(\mB{4},\Z_4) $$
Therefore, with less dependencies than what the wreath product requires, the mod-4 counter can be realized by an isomorphic cascade structure (Fig.\ \ref{fig:mod4counter}).
\end{example}

\begin{example}
The lamplighter group can easily be described as a cascade product:
$$\left(\Z,\Z\right)\wr_T\left(\mathsf{F}(\Z),\bigoplus_{i\in\Z}\Z_2\right)$$
$\mathsf{F}(\Z)$ is the set of all finite subsets of $\Z$. A cascaded state is a 2-tuple $(k,S)$. The top level coordinate is the position of the lamplighter, $k\in\Z$. The second level coordinate is the finite set of lit lamps, $S\in \mathsf{F}(\Z)$. 
The generator set consists of two permutation cascades $T=\{(+0,\ell_2), (+1,t_2)\}$. Lighting the lamps happens on the second level. Which lamp to switch is determined by the position of the lamplighter $\ell=(+0,\ell_2)$, $\ell_2:k\mapsto +1 \text{ of the }k\text{th position}$,
\begin{equation*}
S^{\ell_2(k)} = \begin{cases} S\setminus\{k\}, & \mbox{if } k\in S \\ S\cup\{k\}, & \mbox{if } k\notin S \end{cases}
\end{equation*}
and $t=(+1,t_2)$, where $t_2$ is the constant function giving the identity of the second level group component.
\end{example}

\section{Group Products as Cascade Products}
We show how the usual kinds of products of groups can be described in terms of cascade products.
\subsection{Direct Product}
Using permutation cascades with constant dependency functions we can realize cascade products with no connections at all between the components. Therefore the direct product is a cascade product.

Let $L=\left[(X_1,G_1),\ldots,(X_n,G_n)\right]$.
A \emph{constant dependency function} is defined by $d_{g_i}: (x_1,\ldots,x_{i-1})\mapsto g_i$ for all $(x_1,\ldots,x_{i-1})\in X_1\times\cdots\times X_n$.
For each generator $g_i$ in $G_i$ we define the permutation cascade $$\left(d_{1_{G_1}},d_{1_{G_2}},\ldots,d_{1_{G_{i-1}}},d_{g_i},d_{1_{G_{i+1}}},\ldots,d_{1_{G_n}} \right).$$

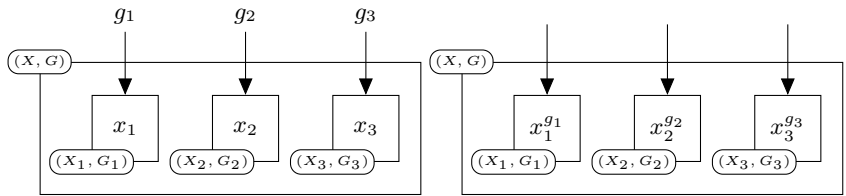
\begin{figure}
\begin{tikzpicture}[node distance=30pt,->,>=triangle 45]
\tikzstyle{component} = [draw=black,rectangle,minimum size=25pt]
\tikzstyle{label} =[draw=black,rectangle, rounded corners,fill=white,inner sep=2pt, inner ysep=2pt]
\node[component] (x1g1) at (0,0) {$x_1$};
\node[label] at (x1g1.south west) {\tiny$(X_1,G_1)$};
\node[component] (x2g2) at (1.6,0) {$x_2$};
\node[label] at (x2g2.south west) {\tiny$(X_2,G_2)$};
\node[component] (x3g3) at (3.2,0) {$x_3$};
\node[label] at (x3g3.south west) {\tiny$(X_3,G_3)$};
\node[component,inner sep=72pt,inner ysep=25pt] at (1.4,0) (xg) {};
\node[label] at (xg.north west) {\tiny$(X,G)$};
\tikzset{node distance=30pt};
\node (s1) at (0,1.5) {$g_1$};
\node (s2) at (1.6,1.5) {$g_2$};
\node (s3) at (3.2,1.5) {$g_3$};
\path[->]  (s1) edge (x1g1);
\path[->]  (s2) edge (x2g2);
\path[->]  (s3) edge (x3g3);
\end{tikzpicture}
\begin{tikzpicture}[node distance=30pt,->,>=triangle 45]
\tikzstyle{component} = [draw=black,rectangle,minimum size=25pt]
\tikzstyle{label} =[draw=black,rectangle, rounded corners,fill=white,inner sep=2pt, inner ysep=2pt]
\node[component] (x1g1) at (0,0) {$x_1^{g_1}$};
\node[label] at (x1g1.south west) {\tiny$(X_1,G_1)$};
\node[component] (x2g2) at (1.6,0) {$x_2^{g_2}$};
\node[label] at (x2g2.south west) {\tiny$(X_2,G_2)$};
\node[component] (x3g3) at (3.2,0) {$x_3^{g_3}$};
\node[label] at (x3g3.south west) {\tiny$(X_3,G_3)$};
\node[component,inner sep=72pt,inner ysep=25pt] at (1.4,0) (xg) {};
\node[label] at (xg.north west) {\tiny$(X,G)$};
\tikzset{node distance=30pt};
\node (s1) at (0,1.5) {$$};
\node (s2) at (1.6,1.5) {$$};
\node (s3) at (3.2,1.5) {$$};
\path[->]  (s1) edge (x1g1);
\path[->]  (s2) edge (x2g2);
\path[->]  (s3) edge (x3g3);
\end{tikzpicture}
\caption{Direct product (parallel composition) of permutation groups. $(X,G)=(X_1,G_1)\times(X_2,G_2)\times(X_3,G_3)=(X_1\times X_2\times X_3, G_1\times G_2\times G_3)$. The coordinatized state $(x_1,x_2,x_3)$ is transformed to $\lp x_1^{g_1},x_2^{g_2},x_3^{g_3}\rp$ under the tuple of permutations $(g_1,g_2,g_3)$, $x_i\in X_i$, $g_i\in G_i$.}
\end{figure}

\begin{figure}
\begin{tikzpicture}[node distance=40pt,->,>=triangle 45]
\tikzstyle{component} = [draw=black,rectangle,minimum size=25pt]
\tikzstyle{label} =[draw=black,rectangle, rounded corners,fill=white,inner sep=2pt, inner ysep=2pt]
\node[fill=black!10] (d1) at(0,4) {$g_1=d_1$};
\node (d2) at (3,4) {$d_2$};
\node (d3) at (6,4) {$d_3$};

\node[fill=black!10] (g2) at (3,2.5) {$g_2=d_2(x_1)$};
\node[fill=black!10] (g3) at (6,1) {$g_3=d_3(x_1,x_2)$};

\node[component] (x1g1) at (0,2.5)  {$x_1$};
\node[label] at (x1g1.south west) {\tiny$(X_1,G_1)$};
\node[component] (x2g2) at (3,1) {$x_2$};
\node[label] at (x2g2.south west) {\tiny$(X_2,G_2)$};

\node[component] (x3g3) at (6,-0.5) {$x_3$};
\node[label] at (x3g3.south west) {\tiny$(X_3,G_3)$};
\node[component,inner sep=130pt,inner ysep=66pt] at (3,1.1) (xg) {};
\node[label] at (xg.north west) {\tiny$(X,G)$};
\tikzset{node distance=30pt};

\path[->]  (d1) edge (x1g1);
\path[>=angle 60,densely dotted]  (x1g1) edge (g2);
\path[->]  (d2) edge (g2);
\path[->]  (g2) edge (x2g2);
\path[->]  (g3) edge (x3g3);
\path[->]  (d3) edge (g3);
\path[>=angle 60,densely dotted]  (x2g2) edge (g3);
\draw[>=angle 60,densely dotted]  (1,2.5) -- (1,0) -- (5,0) -- (5,0.7);

\end{tikzpicture}

\begin{tikzpicture}[node distance=40pt,->,>=triangle 45]
\tikzstyle{component} = [draw=black,rectangle,minimum size=25pt]
\tikzstyle{label} =[draw=black,rectangle, rounded corners,fill=white,inner sep=2pt, inner ysep=2pt]
\node (d1) at(0,4) {$\quad$};
\node (d2) at (3,4) {$$};
\node (d3) at (6,4) {$$};

\node[fill=black!10] (g2) at (3,2.5) {$\quad$};
\node[fill=black!10] (g3) at (6,1) {$\quad\quad\quad\quad$};

\node[component] (x1g1) at (0,2.5)  {$x_1^{g_1}$};
\node[label] at (x1g1.south west) {\tiny$(X_1,G_1)$};
\node[component] (x2g2) at (3,1) {$x_2^{g_2}$};
\node[label] at (x2g2.south west) {\tiny$(X_2,G_2)$};

\node[component] (x3g3) at (6,-0.5) {$x_3^{g_3}$};
\node[label] at (x3g3.south west) {\tiny$(X_3,G_3)$};
\node[component,inner sep=130pt,inner ysep=66pt] at (3,1.1) (xg) {};
\node[label] at (xg.north west) {\tiny$(X,G)$};
\tikzset{node distance=30pt};

\path[->]  (d1) edge (x1g1);
\path[>=angle 60,densely dotted]  (x1g1) edge (g2);
\path[->]  (d2) edge (g2);
\path[->]  (g2) edge (x2g2);
\path[->]  (g3) edge (x3g3);
\path[->]  (d3) edge (g3);
\path[>=angle 60,densely dotted]  (x2g2) edge (g3);
\draw[>=angle 60,densely dotted]  (1,2.5) -- (1,0) -- (5.4,0) -- (5.4,0.9);

\end{tikzpicture}
\caption{Action in a cascade product of components $[(X_1,G_1)$, $(X_2,G_2)$, $(X_3,G_3)]$. The cascaded state $(x_1,x_2,x_3)$ is transformed to $\lp x_1^{ g_1},x_2^{g_2},x_3^{g_3}\rp$ by the permutation cascade $(d_1, d_2, d_3)$. The component actions $g_i$ are calculated by evaluating the dependency functions of $(d_1,d_2,d_3)$ on the states of the components above. The evaluations are highlighted and they happen at the same time. The dependencies, where the state information travels, are denoted by dotted lines.}
\label{fig:cascaction}
\end{figure}
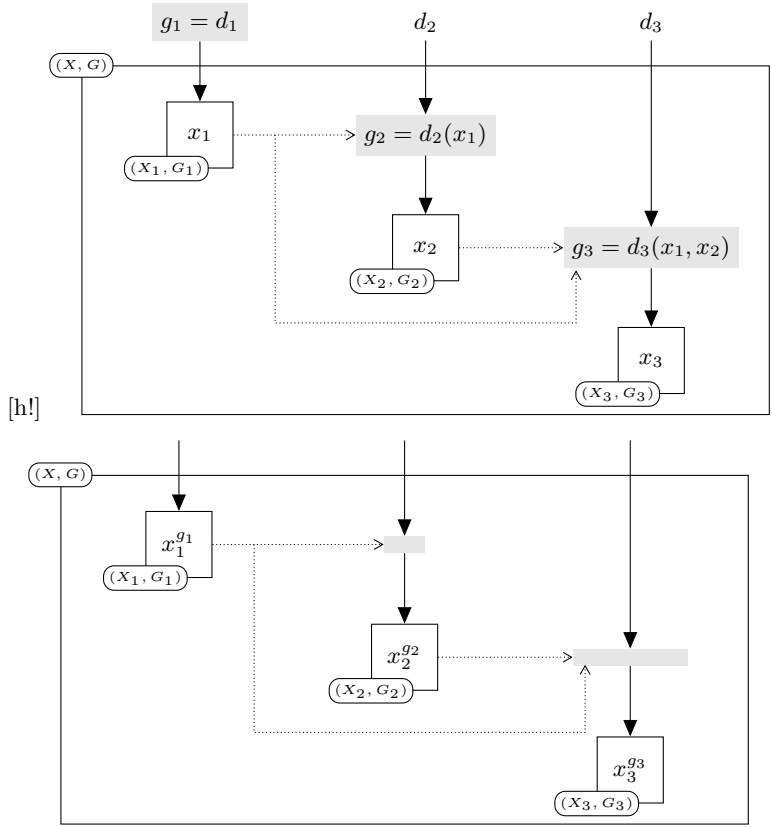

\subsection{Semidirect Product} First we define the semidirect product then show how it can be transformed into a cascade product.
 
\begin{definition}[External Semidirect Product]
Given any groups $H$ and $N$ and a homomorphism $\theta:H\rightarrow\Aut(N)$, denote the automorphism $\theta(h)$ by $\theta_h$. Let $G$ be the set of ordered pairs $\{(h,n)\mid h\in H, n\in N\}$. Then $G$ is a group with the multiplication
$$(h_1,n_1)(h_2,n_2)=(h_1h_2, n_1\theta_{h_1}(n_2))$$
and it is denoted by $H\ltimes_\theta N$.
\end{definition}

The semidirect product is defined for abstract groups, so we have to take the regular representations of the components. The state set simply becomes the set of group elements. Therefore, the underlying list of components is $[(H,H),(N,N)]$. Then we can rewrite the multiplication in the semidirect product into coordinatewise action:
$$(h_1,n_1)^{(h_2,n_2)}=\lp h_1^{h_2}, n_1^{\theta_{h_1}(n_2)}\rp$$
where $\theta_{h_1}(n_2)\in N$ and it can be considered as a function $H\times N\rightarrow N$. For fixed $n_2$,  we can make this into a function $f_2:H\rightarrow N$, which is exactly the bottom level dependency function corresponding to $(h_2,n_2)$.  
Thus, the permutation cascade  $(f_1,f_2)$ with constant function $f_1=h_2$ yields $(h_2,n_2) \in H\ltimes_\theta N$.

\subsection{Wreath Product} The full cascade product is isomorphic to the iterated permutation wreath product. There are two distinct ways to define the wreath product. One way hides, while the other one clearly shows, that we have indeed all dependencies present. 
\label{sec:wreath}

\begin{definition}[Permutation Wreath Product]
Given a permutation group $(X,H)$ and a group $G$, the \emph{permutation wreath product} is the semidirect product $H\ltimes_\theta G^n$, where $n=|X|$ and $\theta:H\rightarrow \text{Aut}(G^n)$ is a homomorphism  mapping $h$ to $\theta_h$. An automorphism $\theta_h$ of $G^n$ is the permutation of the direct factors
$$\theta_h(g_1,g_2,\ldots,g_n)=(g_{(1^h)},g_{(2^h)},\ldots,g_{(n^h)})$$
shuffling the positions according to the action of $h$ on $X$. Therefore it is an action by conjugation on $G^n$ by $H$.
\end{definition}


\begin{theorem}
The permutation wreath product of a permutation group $(X,H)$ and a group $G$ is isomorphic to the group acting in the full cascade product $(X,H)\Bwr(Y,G)$.
\end{theorem}
\PF
Consider functions from set $X$ to group $G$, $|X|=n$. 
Now $G^n$ is the same as all functions from $X$ to $G$ under pointwise multiplication: 
For $f,g:X \rightarrow G$,  $(f\ast g)(x):=f(x)g(x)$ for any $x\in X$, which obviously corresponds exactly to
$$(f(x_1),\ldots, f(x_n)) (g(x_1),\ldots,g(x_n))=(f(x_1)g(x_1),\ldots,f(x_n)g(x_n)) \mbox{ in $G^n$.}$$
So,
$$(G^X,\ast)\cong\underbrace{G\times G\times\cdots\times G}_{|X|\text{ copies of }G}.$$
That is, we can view $f : X \rightarrow G$  as an $n$-tuple $(f(x_1),\cdots, f(x_n))$ in $G^n$,

The permutation cascades of the full cascade product are given by $(h,f)$ where $h\in H$ and  $f: X\rightarrow G$.    We have $(h,f)(h',f') = (hh', f \theta_h(f'))$ 
mapping $(x,y)$ to $\lp x^{hh'}, y ^ {f(x) f'(x^h)}\rp, $ for all $(x,y) \in X\times Y$.   
It follows that $H \ltimes_\theta  G^n$   is the group acting in the full cascade product $(X,H)\Bwr(G,G)$,  or  in fact in $(X,H)\Bwr(Y,G)$  for any permutation group $(Y,G)$.
\QED

\begin{example} What are the cascade products that can be built from $\left[(\mB{2},\Z_2),(\mB{2},\Z_2)\right]$? Now the answer is very simple, the subgroups of $(\mB{2},\Z_2)\wr (\mB{2},\Z_2) \cong (\mB{4},D_4)$.
\end{example}

\begin{example} The quaternion group $Q$ is not a semidirect product, but it embeds into the full cascade product $(\mB{2},\Z_2)\Bwr(\mB{2},\Z_2)\Bwr(\mB{2},\Z_2)$, a group with 128 elements. Therefore, it can be described as a cascade product over the components $[(\mB{2},Z_2),(\mB{2},\Z_2),(\mB{2},\Z_2)]$.

In this case,  the dependency functions can only have two values, thus to define cascade permutations it is enough to give only those arguments that give  $+1$ (the generator of $\Z_2$). A cascade permutation realizing $i$ is defined by the dependency functions $(d_1,d_2,d_3)$ where
$$d_2(0)=d_2(1)=d_3(0,0)=d_3(1,1)=+1$$
and all other arguments map to the identity. 
Similarly, a cascade realizing $j$ is defined by $(d'_1, d'_2,d'_3)$ where $$d'_1(\varnothing)=d'_3(0,0)=d'_3(0,1)=+1.$$
One can check that these two order 4 elements generate the 8-element quaternion group Q. Therefore by $W=\{(d_1,d_2,d_3),(d'_1,d'_2,d'_3)\}$ we have
$$(Q,Q)\cong (\mB{2},\Z_2)\wr_W(\mB{2},\Z_2)\wr_W(\mB{2},\Z_2).$$
Note that conjugating $W$ by any element of the full cascade product yields another cascade product isomorphic to $Q$.
\end{example}

\section{Generalizations and Connections}

Cascade product can be defined for any other action, most notably for transformation semigroups in algebraic automata theory.

Cascade actions can be described as actions on trees, similarly to techniques in infinite group theory \cite{BranchGroups2003,Nekrashevych2005self}. 

The linear order can be generalized to a partial order \cite{GeneralizedWreathProducts1983,NehanivNATO95}.

If the connection network is not unidirectional, then we have a more general theory of automata networks \cite{automatanetworks2005}. Roughly speaking, if feedback connections are allowed it becomes easy to build any structure. For instance group components can be constructed from aperiodic building blocks. 

Cascade product as an operator on collections of groups and semigroup is a combination of the wreath products and taking substructures. See \cite{QBook} for a comprehensive view of other operators.

\subsection*{Acknowledgements}
The current definition of cascade product is meant to be tangible for computer scientists and meaningful for mathematicians. It is the distilled result of many discussions and programming sessions. 
The authors would like to thank  Alfredo Donno and Daniele D'Angeli (for pointing to \cite{GeneralizedWreathProducts1983}), James East (for cleaning up foundational issues at the top level), Murray Elder (for suggesting infinite examples), Andrew Francis (for the figure of dependencies), Roozbeh Hazrat (for suggesting the exponent notation), and James D. Mitchell (for many things coming from his new computational implementation from scratch).
The research leading to these results has received funding from the European Union's Seventh Framework Programme (FP7/2007-2013) under grant agreement no.~318202.
\bibliographystyle{plain}
\bibliography{coords}

\end{document}